# a-T-menability of Baumslag-Solitar groups

Światosław R. Gal
Instytut Matematyczny, Uniwersytet Wrocławski
pl. Grunwaldzki 2/4, 50-384 Wrocław
http://www.math.uni.wroc.pl/~sgal/

Tadeusz Januszkiewicz
Instytut Matematyczny Polskiej Akademii Nauk
and
Instytut Matematyczny, Uniwersytet Wrocławski
pl. Grunwaldzki 2/4, 50-384 Wrocław
http://www.math.uni.wroc.pl/~tjan/

## Abstract

The Baumslag-Solitar groups are a-T-menable. This is proved by embedding them into topological groups and studying representation theoretic properties of the latter.

The paper is motivated by the questions of A. Valette.

**Topological groups and approximation properties.**

We investigate some approximation properties of not necessarily discrete groups. We are mainly interested in Baumslag-Solitar groups which are discrete. However in the process topological groups (Lie groups and the automorphism group of a tree) arise. Therefore in first two sections we clarify some concepts that have been used so far only in context of discrete groups.

**Definition 1 (M.** *Gromov). A locally compact, second countable, compactly generated group $G$ is* a-T-menable *iff there exists a metrically proper affine isometric action of $G$ on some Hilbert space.*

a-T-menability is often referred to as the *Haagerup approximation property*. It is equivalent to existence of an approximate identity consisting of $C_0$ positive definite functions on $G$. For ample discussion one may consult [CCJJV].

2000 *Mathematics Subject Classification*: 20f65

Both authors artially supported by a KBN grant 2 P03A 023 14.

**BS-groups**

Let $G \subset \mathfrak{N}$ be a closed subgroup of a locally compact compactly generated topological group $\mathfrak{N}$. Let $i_k \colon H \to G$ ($k = 1, 2$) be two inclusions onto finite index open subgroups, which are conjugated by an automorphism $\phi$ of $\mathfrak{N}$. The main case of interest is when $G$ is discrete and $\mathfrak{N}$ is Lie group.

**Definition 3.** *The $\mathfrak{N}$-BS group $\Gamma$ is the group derived from $(G, H, i_1, i_2)$ by the (topological) HNN construction. In other words, if $G$ is given by the presentation $<S|R>$, $\Gamma$ has the presentation $<S, t|R, t i_1(g) t^{-1} = i_2(g) \ \forall g \in H>$.*

Since we are working with topological presentations, recall that the topology in the HNN extension $\Gamma$ is given by the prebasis $\mathcal{B} = \{\gamma U \gamma' \colon U \text{ open in } G, \ \gamma, \gamma' \in \Gamma\}$. It is clear, that $\Gamma$ is a topological group w.r.t. $\mathcal{B}$.

To see that $\mathcal{B}$ is in fact a basis consider $U$ open in $G$. We can decompose $U$ with respect to $i_1(H)$ as follows $U = \bigcup i_1(U_n) g_n$, where $g_n$ are representatives of cosets of $i_1(H)$. Then $tU = \bigcup i_2(U_n) t g_n$. Therefore any set of $\mathcal{B}$ can be written uniquely in form $\bigcup U_\gamma \gamma$ where $U_\gamma$ runs over open subsets of $G$ and $\gamma$ runs over the chosen set of representatives of $G \backslash \Gamma$. Thus $\mathcal{B}$ is closed under the intersections, in particular is a basis. The same argument shows that the topology on $G$ induced from $\Gamma$ concides with the oryginal one (ie $G \subset \Gamma$ is an open embedding).

**Examples:**

1. ($\mathfrak{N} = \mathbb{R}$) Baumslag-Solitar group with parameters $p$ and $q$ is given by a following presentation $BS_q^p = \langle x, t | x^p = t x^q t^{-1} \rangle = HNN(\mathbb{Z}, \mathbb{Z}, p \cdot, q \cdot)$.
2. ($\mathfrak{N} = \mathbb{R}^n$) Torsion free, finitely generated, abelian-by-cyclic groups are exactly ascending ($i_1 = id$) HNN extensions of $\mathbb{Z}^n$ with $i_2$ given by a $n$ by $n$ matrix with nonzero determinant [BS].
3. $\mathfrak{N}$ a homogeneous nilpotent group (i.e. one admitting a dilating automorphism $\phi$) with a discrete subgroup $G \subset \mathfrak{N}$ such that $\phi(G) \subset G$.

An obvious adaptation of the Bass-Serre theory [S] to the topological context shows that for a topological HNN extension $\Gamma$ of a group $G$ there is a tree $T$ with an edge-transitive $\Gamma$-action such that the vertex stabilizers are conjugated to $G$ and edge stabilizers are conjugated to $H$.

Condition, that $i_k(H)$ are of finite index in $G$ ensures, that $T$ is locally finite. The simplicial automorphism group $Aut(T)$ carries the natural (compact-open) topology with the basis of neighborhoods of the identity $\mathcal{U}_K := \{g | gv = v \ \forall v \in K\}$ where $K$ runs over the family of the compact subsets $K \subset T$. Denote $j_T \colon \Gamma \to Aut(T)$ the homomorphism given by the action.

Define $\widetilde{\mathfrak{N}} = \mathbb{Z} \ltimes \mathfrak{N}$, to be the semidirect product given by the $\mathbb{Z}$ action on $\mathfrak{N}$ via $\phi$. There is an obvious homomorphism $j_\mathfrak{N} \colon \Gamma \to \widetilde{\mathfrak{N}}$, which is the identity on $G$ and sends $t$ to the generator of $\mathbb{Z}$.

**Theorem 1.** *Let $\Gamma$ be $\mathfrak{N}$-BS group as in Definition 3, and let $j_T$ and $j_\mathfrak{N}$ be the homomorphisms defined above. Then the homomorphism $j = (j_T, j_\mathfrak{N}) \colon \Gamma \to Aut(T) \times \widetilde{\mathfrak{N}}$ is an embedding onto a closed subgroup.*

*Proof* : Observe that since $j_{\mathfrak{N}}$ restricted to $G$ is an embedding onto a closed subgroup (this follows from the fact that $G$ in closed subgroup in $\mathfrak{N}$) the same is true for $j$ restricted to $G$. Moreover, if $v$ is the vertex stabilized by $G$, than $j_T(G) = \text{Stab}_v \cap j_T(\Gamma)$.

Let $j(g) = 1$. Since $j_T(g) = 1$, $gv = v$ ie $g \in G$. Since $j_{\mathfrak{N}}$ restricted to $G$ is embedding, $g = 1$.

Let $\gamma \notin j(\Gamma)$. Since $\Gamma$ acts transitively on $T$, we can multiply $\gamma$ by some element of $j(\Gamma)$ and assume that $\gamma v = v$. Then $(\text{Stab}_v - G) \times \mathfrak{N}$ is an open neighbourhood of $\gamma$ omitting $j(\Gamma)$. So $j(\Gamma)$ is closed.

Topologies on $G$ and $j(G)$ coincide. Since $G$ and $j(G)$ are open in $\Gamma$ and $j(\Gamma)$ respectively, the same is true for $\Gamma$ and $j(\Gamma)$.

**Corollary 1.** *If $\mathfrak{N}$ is a-T-menable then so are $\mathfrak{N}$-BS groups.*

*Proof* : The group $\text{Aut}(T)$ is a-T-menable by a result of Haagerup [Ha].

A cyclic extension of a-T-menable group is a-T-menable by [J]. Therefore $\widetilde{\mathfrak{N}}$ is a-T-menable.

It is clear that a product of a-T-menable groups is also a-T-menable and that a closed subgroup of a-T-menable group is also a-T-menable. Thus the corollary follows from Theorem 1. □

**Remark 1.** According to Chapter 3 of [CCJJV] every a-T-menable connected Lie group $\mathfrak{N}$ is locally isomorphic to a direct product of an amenable group and copies of $SO(1,n)$ and $SU(1,n)$. A Lie group is amenable if it is compact-by-solvable. Thus the examples 1.-3. are a-T-menable.

**Remark 2.** Another proof of Corollary 1, relying on the study of a-T-menable of amalgams, is given in a forthcoming paper by the first author.